\documentclass[11pt]{amsart}
\usepackage{a4wide,amssymb,hyperref}
\usepackage[all]{xy}
\usepackage{color}
\usepackage{todonotes}
\date{}
\newtheorem{theorem}{Theorem}[section]
\newtheorem{lemma}[theorem]{Lemma}
\newtheorem{corollary}[theorem]{Corollary}
\newtheorem{proposition}[theorem]{Proposition}
\newtheorem{question}[theorem]{Question}
\newtheorem{problem}[theorem]{Problem}
\newtheorem{example}[theorem]{Example}
\newtheorem{claim}[theorem]{Claim}
\theoremstyle{definition}
\newtheorem{definition}[theorem]{Definition}

\def\R{\mathbb R}
\def\Z{\mathbb Z}

\def\T{\mathbb T}
\def\K{\mathcal K}

\def\U{\mathcal U}
\def\ol{\overline}
\newcommand{\w}{\omega}
\newcommand{\e}{\varepsilon}
\newcommand{\IR}{\mathbb R}
\newcommand{\IZ}{\mathbb Z}
\newcommand{\IN}{\mathbb N}
\newcommand{\IF}{\mathbb F}
\newcommand{\IT}{\mathbb T}
\newcommand{\IH}{\mathbb H}
\newcommand{\IC}{\mathbb C}
\newcommand{\pr}{\mathrm{pr}}
\newcommand{\Ra}{\Rightarrow}

\begin{document}

\title[On unconditionally convergent series in topological rings]
{On unconditionally convergent series\\ in  topological rings}
\author[Taras Banakh]{Taras Banakh}
\address{Taras Banakh: Ivan Franko National University of Lviv (Ukraine), and Jan Kochanowski University in Kielce (Poland)}
\email{t.o.banakh@gmail.com}
\author[Alex Ravsky]{Alex Ravsky}
\address{Alex Ravsky: Pidstryhach Institute for Applied Problems of Mechanics and Mathematics
National Academy of Sciences of Ukraine}
\email{alexander.ravsky@uni-wuerzburg.de}
\keywords{topological ring, unconditional convergence, locally compact topological ring,
locally compact Abelian topological group}
\subjclass{13J10, 22A10, 54H11, 54H13}
\begin{abstract}
We define a topological ring $R$ to be \emph{Hirsch}, if for any unconditionally convergent series $\sum_{n\in\w} x_i$ in $R$ and any
neighborhood $U$ of the additive identity $0$ of $R$ there exists a neighborhood $V\subseteq R$ of
$0$ such that $\sum_{n\in F} a_n x_n\in U$ for any finite set $F\subset\w$ and any sequence
$(a_n)_{n\in F}\in V^F$. We recognize Hirsch rings in certain known classes of
topological rings. For this purpose we introduce and develop the technique of seminorms on actogroups.
We prove, in particular, that a topological ring $R$ is Hirsch provided
$R$ is locally compact or $R$ has a base at the zero consisting of
open ideals or $R$ is a closed subring of the Banach ring $C(K)$, where $K$ is a compact Hausdorff
space. This implies that the Banach ring $\ell_\infty$ and its subrings $c_0$ and $c$ are Hirsch. Also we prove that for every $p\in[1,2]$ the Banach ring $\ell_p$ is Hirsch. On the other hand, for any distinct numbers $p,q\in[1,\infty]$ the commutative Banach ring $\ell_p\oplus i\ell_q$ is not Hirsch. Also for any $p\in (1,\infty)$, the (noncommutative) Banach ring $L(\ell_p)$ of continuous
endomorphisms of the Banach ring $\ell_p$ is not Hirsch. We do not know whether the Banach rings $\ell_p$ are Hirsch for $p\in(2,\infty)$.
\end{abstract}

\maketitle \baselineskip15pt

\section{Introduction}


We define a topological ring $R$ to be \emph{Hirsch} if for any unconditionally convergent series
$\sum_{n\in\w} x_n$ in $R$ and any neighborhood $U\subseteq R$ of the additive identity $0$ of $R$
there exists a neighborhood $V\subseteq R$ of $0$ such that $\sum_{n\in F} a_nx_n \in U$ for any
finite set $F\subset \w$ and any sequence $(a_n)_{n\in F}\in V^F$.


This notion was suggested by a question~\cite{Dch} of Dylan Hirsch  who was motivated by the idea of generalizing the standard theory of Lebesgue measure and integral to functions with values in topological rings.

 In this paper we consider the following general problem.


\begin{problem}[see~\cite{Dch}]\label{prob} Recognize Hirsch rings in known classes of topological rings.
\end{problem}


 Exploring Problem~\ref{prob}, we shall introduce and develop the technique of seminorms on actogroups, see Section~\ref{sec:HA}. We
shall show, in particular, that a topological ring $R$ is Hirsch provided   $R$ has a base at the zero consisting of open ideals (see Proposition~\ref{p4}) or $R$ is locally compact (see Theorem~\ref{t:lc-Hr})
 or $R$ is a subring of the topological rings $C_{\mathcal K}(X,S)$ of continuous functions from a topological space $X$ to a locally compact topological ring $S$, endowed with the topology of uniform convergence on sets from an ideal $\mathcal K$ of compact Hausdporff subsets of $X$ (see Theorem~\ref{t:CK-ring}).
  This implies that the Banach ring $\ell_\infty$ and its subrings $c_0$ and $c$ are Hirsch, see Theorem~\ref{t:6.3}.     Also we prove that for every $p\in[1,2]$ the Banach ring $\ell_p$ is Hirsch, see  Corollary~\ref{c:12}. On the other hand, for any distinct numbers $p,q\in[1,\infty]$ the commutative Banach ring $\ell_p\oplus i\ell_q$ is not Hirsch, see Theorem~\ref{t:pq}. Also for any $p\in (1,\infty)$, the (noncommutative) Banach ring $L(\ell_p)$ of continuous
endomorphisms of the Banach ring $\ell_p$ is not Hirsch, see Propositions~\ref{p1} and~\ref{p1'}. We do not know whether the Banach rings $\ell_p$ are Hirsch for $p\in(2,\infty)$, see Problem~\ref{prob2}.

\section{Unconditional convergence in topological groups}

A {\em group} in this paper means an Abelian topological group. We shall use the additive notation for denoting the group operations. For a group $G$ its identity element is denoted by $0_G$ (so, $x+0_G=x=0_G+x$ for every $x\in G$).

Let $I$ be a countable infinite set and $(x_i)_{i\in I }$ be a sequence of elements of a
 group $G$. The sequence $(x_i)_{i\in I}$
\begin{itemize}
\item  \emph{converges} to an element $x\in G$
if for any neighborhood $U$ of the identity of $G$ there exists a finite set $F\subset I$ such that
$x_i\in x+U$ for any $i\in I\setminus F$;
\item is \emph{Cauchy}
if for any neighborhood $U$ of the identity of $G$ there exists a finite set $F\subset I$ such that
$x_i-x_j\in U$ for any $i,j\in I\setminus F$.
\end {itemize}
It is easy to show that any convergent sequence is Cauchy.
A group $G$ is {\it sequentially complete} if any Cauchy sequence in $G$ converges.

The series $\sum_{i\in I} x_i$
\begin{itemize}
\item  \emph{unconditionally converges} to an element $x\in G$ if for
any neighborhood $U$ of the identity of $G$ there exists a finite set $F\subset I$ such that
$\sum_{i\in J}x_i\in x+U$ for any finite set $J\subset I$ containing $F$
(if the set $J$ is empty then we put $\sum_{i\in J}x_i=0$);
\item is \emph{unconditionally Cauchy} if for
any neighborhood $U$ of the identity of $G$ there exists a finite set $F\subset I$ such that
$\sum_{i\in K}x_i\in U$ for any finite set $K\subset I\setminus F$.
\end{itemize}
It is easy to show that any unconditionally convergent series is unconditionally Cauchy.
On the other hand we have the following

\begin{proposition}\label{ConToUnCon} Let $\sum_{i\in I} x_i$ be an unconditionally Cauchy
series in a group $G$ and  $(I_n)_{n\in\w}$ be a nondec\-reasing sequence of finite sets such
that $\bigcup_{n\in\w} I_n=I$. If the sequence $\left(\sum_{i\in I_n}x_i\right)_{n\in\w}$ converges
to a point $x\in G$ then the series $\sum_{i\in I} x_i$ unconditionally converges to $x$.
\end{proposition}

\begin{proof} Let $U$ be any neighborhood of the identity of $G$. Pick a neighborhood $V$ of
the identity of $G$ such that $V+V\subseteq U$.
Since the sequence $\left(\sum_{i\in I_n}x_i\right)_{n\in\w}$ converges to $x$,
there exists $m\in\w$, such that $\sum_{i\in I_n}x_i\in x+V$ for any $n\ge m$.
Since the series $\sum_{i\in I} x_i$ is
unconditionally Cauchy, there exists a finite set $J\subset I$ such that
$\sum_{i\in K}x_i\in V$ for any finite set $K\subset I\setminus J$.
Since the sequence $(I_n)_{n\in\w}$ is nondecreasing and $\bigcup_{n\in\w} I_n=I$,
there exists a number $k\ge m$ such that $I_k\supseteq J$.
Let $F\subset I$ be any finite set containing $I_k$. Then
$$\sum_{i\in F} x_i=\sum_{i\in I_k} x_i+\sum_{i\in F\setminus I_k} x_i\in x+V+V\subseteq x+U.$$
\end{proof}

\begin{proposition}\label{CauchyToUnCon} Any unconditionally Cauchy series $\sum_{i\in I} x_i$ in
a sequentially complete group $G$ is unconditionally convergent.
\end{proposition}

\begin{proof} Fix any increasing sequence $(I_n)_{n\in\w}$ of finite subsets of
$I$ such that $\bigcup_{n\in\w} I_n=I$. Let $U$ be any neighborhood of the identity of $G$
such that $U=-U$. Since the series $\sum_{i\in I} x_i$ is unconditionally Cauchy, there exists
a finite set $F\subset I$ such that $\sum_{i\in K}x_i\in U$ for any finite set
$K\subset I\setminus F$.
Since the sequence $(I_n)_{n\in\w}$ is increasing and $\bigcup_{n\in\w} I_n=I$,
there exists $k\in\w$ such that $I_k\supseteq F$. Let $m,n\in\w$ be any numbers with $k\le m\le n$. Since $I_n\setminus I_m\subset I\setminus F$, we have
$$\sum_{i\in I_n} x_i-\sum_{i\in I_m} x_i= \sum_{i\in I_n\setminus I_m} x_i\in U.$$
Thus the sequence $\left(\sum_{i\in I_n}x_i\right)_{n\in\w}$ is Cauchy.
Since the group $G$ is sequentially complete, the sequence $\left(\sum_{i\in I_n}x_i\right)_{n\in\w}$
converges to a point $a\in G$. By Proposition~\ref{ConToUnCon}, the series $\sum_{i\in I} x_i$
unconditionally converges to $x$.
\end{proof}

An {\em endomorphism} of a group $G$ is a continuous homomorphism from $G$ to itself.

\begin{lemma}\label{lConvNew1}
Let $\sum_{i\in I} x_i$ be an unconditionally Cauchy series in a group $G$
and $(v_i)_{i\in I}$ be a sequence of endomorphisms of $G$ such that
the set $D=\{v_i\}_{i\in I}$  is finite.
Then the series $\sum_{i\in I} v_i(x_i)$ is unconditionally Cauchy.
\end{lemma}
\begin{proof} Let $U$ be any neighborhood of the identity of $G$. Pick a neighborhood $V$ of
the identity of $G$ such that the sum of $|D|$ summands $V$ is contained in $U$.
Let $v\in D$ be any element. Put $I_v=\{i\in I: v_i=v\}$ and
pick a neighborhood $W_v$ of the identity of $G$ such that $v(W_v)\subseteq V$.
Put $W=\bigcap_{v\in D} W_v$.
Since the series $\sum_{i\in I} x_i$ is unconditionally Cauchy, there exists
a finite set $F\subset I$ such that $\sum_{i\in K}x_i\in W$ for any finite set
$K\subset I\setminus F$. Then for any finite set $K\subset I\setminus F$ we have
$$\sum_{i\in K}v_i(x_i)=\sum_{v\in D}\sum_{i\in K\cap I_v} v(x_i)=
\sum_{v\in D}v\Big(\sum_{i\in K\cap I_v} x_i\Big)\in \sum_{v\in D}v(W)\subseteq
\sum_{v\in D}v(W_v)\subseteq \sum_{v\in D}V\subseteq U,$$
witnessing that the  series $\sum_{i\in I} v_i(x_i)$ is unconditionally Cauchy.
\end{proof}

\section{Topological rings that have a local base of open ideals}

A {\em topological ring} is a group $R$ endowed with a continuous associative binary
operation $\cdot:R\times R\to R$ which is {\em distributive}, that is
$$x(y+z)=xy+xz\quad\mbox{and}\quad (x+y)z=xz+yz $$ for all $x,y,z\in R$.
 A nonempty subset $I$ of a ring $R$ is called {\em an ideal} in $R$ if $IR\cup RI\subseteq R$.

\begin{proposition}\label{p4} Every topological ring $R$ that has a neighborhood base at zero consisting of open
ideals is Hirsch.
\end{proposition}
\begin{proof} Let $\sum_{i\in\w} x_i$ be any (unconditionally) convergent series in $R$ and $U$ be any neighborhood
of the zero $0$ of $R$. Let $V$ be an arbitrary open ideal contained in $U$. Then for every finite
set $I\subset\w$ we have $\sum_{i\in I}Vx_i\subseteq V\subseteq U$ as $V$ is an ideal in $R$.
\end{proof}

\section{Hirsch actogroups}\label{sec:HA}

It is natural to explore Problem~\ref{prob} in a wider context of active groups, i.e., groups endowed with an action of a topological space with  zero.

\begin{definition}\label{d:S} An {\em  active group} (briefly an {\em actogroup}) $A\curvearrowright G$ is a topological group $G$ endowed with a continuous action $A\times G\to G$, $(a,x)\mapsto ax$, of a topological space $A$ with a distinguished point $\mathsf 0_A\in A$. The action should satisfy two axioms:
\begin{itemize}
\item[(i)] $\mathsf 0_A\cdot x=0_G$ for every $x\in G$;
\item[(ii)] $a(x+x)=ax+ax$ for every $a\in A$ and $x\in G$.
\end{itemize}
The topological space $A$ is called the {\em acting space} of the actogroup $A\curvearrowright G$ and the distinguished point $\mathsf 0_A$ is called {\em the zero} of $A$. If the acting space $A$ and the action $A\times G\to G$ are clear from the context, then we shall write $G$ instead of $A\curvearrowright G$.
\end{definition}

\begin{example} Each topological ring $R$ is an actogroup $R\curvearrowright R$ with respect to the continuous action $R\times R\to R$, $(a,x)\mapsto ax$, assigning to a pair of elements $a,x\in R$ their product $ax$ in the ring $R$.
\end{example}

\begin{definition} An actogroup $A\curvearrowright G$ is defined to be {\em Hirsch} if for every unconditionally convergent series $\sum_{i\in \w}x_i$ in $G$ and every neighborhood $U\subseteq G$ of $0_G$ there exists a neighborhood $V$ of zero $\mathsf 0_A$ in the acting space $A$ of $G$ such that for every finite set $F\subseteq \w$ and sequence $(a_i)_{i\in F}\in V^F$ we have $\sum_{i\in F}a_ix_i\in U$.
\end{definition}

It follows that a topological ring is Hirsch if and only if it is Hirsch as an actogroup. So, Problem~\ref{prob} is a partial case of the following more general

\begin{problem}\label{prob:gen} Recognize Hirsch actogroups in known classes of actogroups.
\end{problem}

In the following subsections we shall present examples of Hirsch actogroups among actogroups whose topology is generated by some specific seminorms.

\subsection{Seminorms on groups}

\begin{definition} A {\em seminorm} on a (topological) group $G$ is a (continuous) function $\|\cdot\|:G\to[0,\infty)$ such that
$$\|0_G\|=0,\quad \|-x\|=\|x\|\quad\mbox{and}\quad \|x+y\|\le\|x\|+\|y\|$$
for every $x,y\in G$.

A seminorm $\|\cdot\|:G\to[0,\infty)$ is a {\em norm} if
\begin{itemize}
\item $\|x\|>0$ for any $x\in G\setminus\{0_G\}$ and
\item $\|nx\|=n\|x\|$ for any $x\in G$ and $n\in\w$.
\end{itemize}
\end{definition}
Here $nx$ is the element of $G$, defined by the recursive formula:
$$\mbox{$0x=0_G$ \ and \ $(n+1)x=nx+x$ for every $n\in\w$.}$$

We say that the topology of a group $G$ is {\em generated by a family} $\mathcal S$ of seminorms if for every neighborhood $U\subseteq G$ of $0_G$ there exists a seminorm $\|\cdot\|\in\mathcal S$ such that $\{x\in G:\|x\|<\e\}\subseteq U$ for some $\e>0$ (by the definition, all seminorms on a  group are  continuous). If the family $\mathcal S$ consists of a single seminorm $\|\cdot\|$, then we say that the topology of $G$ is {\em generated by the seminorm} $\|\cdot\|$.

A seminorm $\|\cdot\|:G\to [0,\infty)$ is {\em bounded} if $\sup_{x\in G}\|x\|<\infty$.
Theorem 3.3.9 from~\cite{ArhTka} implies that the topology of any group is generated by a suitable family of bounded seminorms.

\subsection{Locally $2$-homogeneous seminorms on groups}

\begin{definition} A seminorm $\|\cdot\|$ on a group $G$ is called
\begin{itemize}
\item {\em $2$-homogeneous} if $\|2x\|=2\|x\|$ for every $x\in G$;
\item {\em locally $2$-homogeneous} if there exists $\e>0$ such that $\|2x\|=2\|x\|$ for every $x\in G$ with $\|x\|\le\e$.
\end{itemize}
\end{definition}
Each 2-homogeneous seminorm on a group is locally 2-homogeneous.

\begin{example}\label{ex:R2} Each norm on a group is $2$-homogeneous. In particular, the norm $\|\cdot\|:\IR\to[0,\infty)$, $\|\cdot\|:x\mapsto |x|$, on the additive group of real numbers $\IR$ is $2$-homogeneous.
\end{example}

Next we explore some operations over (locally) $2$-homogeneous seminorms. The following two lemmas can be easily derived from the definitions.

\begin{lemma}\label{lem:hom2hom} Let $h:H\to G$ be a continuous homomorphism between groups. If $\|\cdot\|_G$ is a (locally) $2$-homogeneous seminorm on $G$, then the seminorm
$$\|\cdot\|_H:H\to[0,\infty),\quad \|\cdot\|_H:x\mapsto \|h(x)\|_G,$$ on the group $H$ is (locally) $2$-homogeneous.
\end{lemma}


We recall that for a group $G$ and a subgroup $H\subseteq G$ the quotient group $G/H$ is the group whose elements are cosets $x+H$ where $x\in G$.

\begin{lemma}\label{lem:prod2hom} Let $\|\cdot\|_G$ and $\|\cdot\|_H$ be seminorms on groups $G$ and $H$, respectively. If the seminorms $\|\cdot\|_G$ and $\|\cdot\|_H$ are (locally) $2$-homogeneous, then so is the seminorm $\|\cdot\|:G\times H\to[0,\infty)$, $\|\cdot\|:(x,y)\mapsto\max\{\|x\|_G,\|y\|_H\}$.
\end{lemma}


\begin{lemma}\label{l:h2h} Let $\|\cdot\|_G$ be a locally $2$-homogeneous seminorm on a group $G$ and $H$ be a  subgroup of $G$ such that $\{x\in H:\|x\|_G< 5\e\}=\{0_H\}$ for some $\e>0$. Then the seminorm
$$\|\cdot\|:G/H\to[0,\infty),\quad \|y\|=\inf\{\|x\|_G:x\in y\}$$ is locally $2$-homogeneous.
\end{lemma}

\begin{proof} Replacing $\e>0$ by a smaller positive number, we can assume that $\|2x\|_G=2\|x\|_G$ for all $x\in G$ with $\|x\|_G<\frac32\e$. We claim that for every $y\in G/H$ with $\|y\|\le\e$ we have $\|2y\|=2\|y\|$. Indeed, assume for a contradiction that $\|2y\|\ne 2\|y\|$.  Since $\|2y\|\le \|y\|+\|y\|=2\|y\|$, we have $\|2y\|<2\|y\|$, so there exists $z\in 2y$ such that $\|z\|_G<2\|y\|\le 2\e$.
Since $\|y\|\le \e$, there exists $x\in y$ such that $\|x\|_G<\frac32\e$ and then $\|2x\|_G=2\|x\|_G<3\e$.  It follows from $z,2x\in 2y$ that $z-2x\in H$. The triangle inequality implies that
$$\|2x-z\|_G\le\|2x\|_G+\|z\|_G<3\e+2\e=5\e
$$
and hence $2x-z\in\{h\in H:\|h\|<5\e\}=\{0_H\}$ and finally $z=2x$. Then $$2\|y\|\le 2\|x\|_G=\|2x\|_G=\|z\|_G<2\|y\|,$$ which is a required contradiction.
\end{proof}

\begin{corollary}\label{cor:T2hom} The seminorm $\|\cdot\|:\IR/\IZ\to[0,\infty)$, $\|\cdot\|:y\mapsto\min\{|x|:x\in y\}$, on the quotient group $\IT=\IR/\IZ$ is locally $2$-homogeneous.
\end{corollary}

\begin{lemma}\label{l:1f} Let $H$ be an open subgroup of a group $G$. If a seminorm $\|\cdot\|_H:H\to[0,\infty)$ on the group $H$ is locally $2$-homogeneous, then the seminorm $$\|\cdot\|_G:G\to[0,\infty),\quad \|\cdot\|_G\mapsto\begin{cases}\min\{1,\|x\|_H\}&\mbox{if $x\in H$};\\
1&\mbox{if $x\in G\setminus H$};
\end{cases}
$$on the group $G$ is locally $2$-homogeneous, too.
\end{lemma}

\begin{proof} Since $\|\cdot\|_H$ is locally $2$-homogeneous, there exists a positive $\e<\frac12$ such that $\|2x\|_H=2\|x\|_H$ for every $x\in H$ with $\|x\|_H\le \e$. The definition of the norm $\|\cdot\|_G$ ensures that for every $x\in G$ with $\|x\|_G\le\e$ we will have $x\in H$ and hence $1>\e\ge \|x\|_G=\min\{1,\|x\|_H\}$ implies that $\|x\|_H=\|x\|_G$ and $\|2x\|_H\le 2\|x\|_H\le 2\e<1$. Then
$$\|2x\|_G=\min\{1,\|2x\|_H\}=\|2x\|_H=2\|x\|_H=\min\{2,2\|x\|_H\}=2\min\{1,\|x\|_H\}=2\|x\|_G,$$
so the seminorm $\|\cdot\|_G$ is locally $2$-homogeneous.
\end{proof}


\subsection{Lipschitz and contracting seminorms on actogroups}

By a {\em seminorm} on an actogroup $A\curvearrowright G$ we understand a seminorm on the group $G$.

\begin{definition} A seminorm $\|\cdot\|:G\to[0,\infty)$ on an actogroup $A\curvearrowright G$  is called
\begin{itemize}
\item {\em Lipschitz} if there exist positive real numbers $L,\e$ and a neighborhood $V\subseteq A$ of $\mathsf 0_A$ such that $\|ax\|\le L \|x\|$ for every $a\in V$ and every $x\in G$ with $\|x\|\le\e$;
\item {\em contracting} if for every $\lambda>0$ there exist $\e>0$ and a neighborhood $V\subseteq A$ of $\mathsf 0_A$ such that $\|ax\|\le\lambda \|x\|$ for every $a\in V$ and every $x\in G$ with $\|x\|\le\e$.
\end{itemize}
\end{definition}

It is clear that each contracting seminorm on an actogroup is Lipschitz. 

\begin{theorem}\label{t:Ravsky} Let $\|\cdot\|:G\to[0,\infty)$ be a locally $2$-homogeneous seminorm on an actogroup $A\curvearrowright G$
and assume that there exists $\e>0$ such that the closed $\e$-ball $B=\{x\in G:\|x\|\le\e\}$ is compact. Then  the seminorm $\|\cdot\|$ is contracting.
\end{theorem}
\begin{proof} 
Since $\|\cdot\|$ is locally $2$-homogeneous, we can replace $\e$ by a smaller positive number and assume that  $\e\le 1$ and $\|2x\|=2\|x\|$ for every $x\in B$. To prove that $\|\cdot\|$ is contracting, take any positive real number $\lambda\le 1$.

 The continuity of the action $A\times G\to G$, $(a,x)\mapsto ax$, implies that $W=\{(a,x)\in A\times B:\|ax\|<\tfrac14\lambda\e\}$ is an open neighborhood of the compact set $\{\mathsf 0_A\}\times B$ in $A\times B$. Using the compactness of $B$, we can find a neighborhood $V\subseteq A$ of $0_A$ such that $V\times B\subseteq W$ and hence $\|ax\|<\tfrac14\lambda\e$ for every $a\in V$ and $x\in B$.

To see that $\|\cdot\|$ is contracting, it remains to show that $\|ax\|\le\lambda\|x\|$ for every $a\in V$ and $x\in B$. Suppose for a contradiction that there exist $a\in V$ and $x\in B$ such that  $\|ax\|>\lambda\|x\|$.

It follows from $x\in B$ that $\|x\|\le\e$ and hence $0\le\lambda\|x\|<\|ax\|<\tfrac14\lambda\e\le\e$. Then there exists a unique number $k\in\w$ such that $2^k\|ax\|\le \e<2^{k+1}\|ax\|$.

\begin{claim}\label{cl1} For every nonnegative integer $i\le k+1$ we have $\|a(2^ix)\|=2^i\|ax\|$.
\end{claim}

\begin{proof} We prove the claim by induction on $i$. For $i=0$ the equality $\|a(2^0x)\|=2^0\|ax\|$ is trivial. Assume that $\|a(2^ix)\|=2^i\|ax\|$ for some $i<k+1$. Then $$\|a(2^ix)\|=2^i\|ax\|\le 2^k\|ax\|\le\e$$ and hence $$\|a(2^{i+1}x)\|=\|a(2^ix+2^ix)\|=\|a(2^ix)+a(2^ix)\|=2\|a(2^ix))\|=2^{i+1}\|ax\|$$by the inductive hypothesis.
\end{proof}

\begin{claim} $\|x\|>0$.
\end{claim}

\begin{proof} To derive a contradiction, assume that $\|x\|=0$. Claim~\ref{cl1} and the choice of $k$ imply $\|a(2^{k+1}x)\|=2^{k+1}\|ax\|>\e$. On the other hand, $\|x\|=0$ implies $\|2^{k+1}x\|=0$ and hence $2^{k+1}x\in B$ and  $\|a(2^{k+1}x)\|<\tfrac14\lambda\e\le\e$ by the choice of the neighborhood $V\ni a$. This contradiction shows that $\|x\|>0$.
\end{proof}

Since $0<\|x\|\le\e$, there exists a unique  number $n\in\w$ such that $2^n\|x\|\le\e<2^{n+1}\|x\|$. By analogy with Claim~\ref{cl1} we can prove that $\|2^{i}x\|=2^{i}\|x\|$ for all nonnegative integers $i\le n+1$.

Since $\lambda\le 1$, there exists a unique number $d\in\w$ such that $2^{d}\lambda\le1<2^{d+1}\lambda$.

Observe that for every $i\le n$ we have $\|2^ix\|\le2^i\|x\|\le 2^n\|x\|\le\e$ and hence $2^ix\in B$ and $\|a(2^ix)\|<\tfrac14\lambda\e\le\e<2^{k+1}\|ax\|=\|a(2^{k+1}x)\|$ by the choice of the neighborhood $V\ni a$ and the number $k$. Also for every integer
$i\le d+2$, by the triangle inequality, $\|a(2^{n+i}x)\|\le 2^{i}\|a(2^nx)\|<2^{d+2}\tfrac14\lambda\e\le\e<2^{k+1}\|ax\|=\|a(2^{k+1}x)\|$. Therefore, for every $i\le n+d+2$ we have $\|a(2^ix)\|\le \e<\|a(2^{k+1}x)|$, which implies that $n+d+2<k+1$.
Then Claim~\ref{cl1} and the choice of $k$ and $n$ ensure that
$$2^{n+d+2}\|ax\|=\|a(2^{n+d+2}x)\|\le \e<2^{n+1}\|x\|$$ and finally $\|ax\|<\frac1{2^{d+1}}\|x\|<\lambda\|x\|$, which contradicts our assumption.
\end{proof}

\subsection{Absolute seminorms}

\begin{definition} A seminorm $\|\cdot\|$ on a group $G$ is defined to be
\begin{itemize}
\item  {\em $C$-absolute} for some positive real number $C$ if there exists $\e>0$ such that for any finite set $F$ and sequence $(x_i)_{i\in F}\in G^F$ with $M:=\max_{E\subseteq F}\|\sum_{i\in E}x_i\|\le \e$ we have $\sum_{i\in F}\|x_i\|\le CM$;
\item {\em absolute} if $\|\cdot\|$ is $C$-absolute for some positive  real number $C$.
\end{itemize}
\end{definition}

\begin{example}\label{ex:R} The norm $|\cdot|:\IR\to[0,\infty)$, $|\cdot|:x\mapsto |x|$, on the additive group $\IR$ of real numbers is $2$-absolute.
\end{example}

\begin{proof} Given any finite set $F$ and sequence $(x_i)_{i\in F}\in\IR^F$, consider the sets $F_+=\{i\in F:x_i\ge 0\}$ and $F_-=F\setminus F_+$. It follows from $$\sum_{i\in F}|x_i|=\sum_{i\in F_+}x_i-\sum_{i\in F_-}x_i=\big|\sum_{i\in F_+}x_i\big|+\big|\sum_{i\in F_-}x_i\big|$$ that $$\sum_{i\in F}|x_i|\le2\max\big\{\big|\sum_{i\in F_+}x_i\big|,\big|\sum_{i\in F_-}x_i\big|\big\}\le 2\max_{E\subseteq F}\big|\sum_{i\in E}x_i\big|,$$witnessing that the norm $|\cdot|$ on $\IR$ is 2-absolute.
\end{proof}

A result of Netuka and Vesel\'y \cite{NV} implies the following general fact.

\begin{example}\label{ex:NV} For every positive integer $n$, the standard Euclidean norm $\|\cdot\|:\IR^n\to[0,\infty)$ on the additive group $\IR^n$ is $C$-absolute for the constant
$$C=2\sqrt{\pi}\,\frac{\Gamma(\tfrac{n+1}2)}{\Gamma(\tfrac{n}2)},$$which is the best possible.
\end{example}

As partial cases of Example~\ref{ex:NV} we have two next examples.

\begin{example} The norm $\|\cdot\|:\IC\to[0,\infty)$, $\|\cdot\|:z\mapsto |z|$, on the additive group of complex numbers $\IC$ is $\pi$-absolute.
\end{example}

\begin{example} The  norm $\|\cdot\|:\IH\to[0,\infty)$, $\|\cdot\|:q\mapsto |q|$, on the additive group of quaternion numbers $\IH$ is $3\pi$-absolute.
\end{example}

\begin{lemma}\label{l:prod-a} Let $G,H$ be two groups and  $\|\cdot\|_G,\|\cdot\|_H$ be seminorms on $G,H$, respectively. If for some positive real number $C$ these seminorms  are $C$-absolute, then the norm
$$\|\cdot\|:G\times H\to[0,\infty),\quad \|\cdot\|:(g,h)\mapsto\max\{\|g\|_G,\|h\|_H\},$$
is $2C$-absolute.
\end{lemma}

\begin{proof} Let $F$ be a finite set and $(z_i)_{i\in F}$ be a sequence in $G\times H$. Write each $z_i$ as $(x_i,y_i)$ for some $x_i\in G$ and $y_i\in H$. For the sets $I=\{i\in F:\|x_i\|_G\le\|y_i\|_H\}$ and $J=\{i\in F:\|x_i\|_G>\|y_i\|_H\}$ we have $$
\sum_{i\in F}\|z_i\|=\sum_{i\in F}\max\{\|x_i\|_G,\|y_i\|_H\}=\sum_{i\in J}\|x_i\|_G+\sum_{i\in I}\|y_i\|_H$$ and hence
$$\max\big\{\sum_{i\in J}\|x_i\|_G,\sum_{i\in I}\|y_i\|_H\big\}\ge\frac12\sum_{i\in F}\|z_i\|.$$ Since the norms $\|\cdot\|_G,\|\cdot\|_H$ are $C$-absolute, there are subsets $I'\subseteq  I$ and $J'\subseteq J$ such that $$\big\|\sum_{i\in I'}x_i\big\|_G\ge \frac1C\sum_{i\in  I}\|x_i\|_G\quad\mbox{and}\quad\big\|\sum_{i\in J'}y_i\big\|_H\ge\frac1C\sum_{i\in J}\|y_i\|_H.$$
Then $$\max_{E\subseteq F}\big\|\sum_{i\in E}z_i\big\|\ge\max\big\{\big\|\sum_{i\in I'}x_i\big\|_G,\big\|\sum_{i\in J'}y_i\big\|_H\big\}\ge \frac1C\max\big\{\sum_{i\in I}\|x_i\|_G,\sum_{i\in J}\|y_i\|_H\big\}\ge\frac1{2C}\sum_{i\in F}\|z_i\|,$$witnessing that the norm $\|\cdot\|$ on $G\times H$ is $2C$-absolute.
\end{proof}

Lemma~\ref{l:prod-a} implies

\begin{corollary}\label{c:prod-a} Let $G,H$ be two groups and  $\|\cdot\|_G,\|\cdot\|_H$ be seminorms on $G,H$, respectively. If the seminorms $\|\cdot\|_G ,\|\cdot\|_H$ are absolute, then so is the seminorm
$$\|\cdot\|:G\times H\to[0,\infty),\quad \|\cdot\|:(g,h)\mapsto\max\{\|g\|_G,\|h\|_H\}.$$
\end{corollary}

\begin{lemma}\label{p:quot} Let $C$ be a positive real number, $\|\cdot\|_G:G\to [0,\infty)$ be a $C$-absolute seminorm on a group $G$, and $H$ be a subgroup of $G$ such that $\{x\in H:\|x\|_G<5\e\}=\{0_H\}$ for some $\e>0$. Then the seminorm $$\|\cdot\|:G/H\to[0,\infty),\quad \|\cdot\|:y\mapsto\inf\{\|x\|_G:x\in y\},$$
is $C$-absolute.
\end{lemma}

\begin{proof} Replacing the number $\e$ by a smaller number, if necessary, we can assume that for any finite set $F$ and any sequence $(x_i)_{i\in F}\in G^F$ with $M:=\max_{A\subseteq F}\big\|\sum_{i\in A}x_i\big\|_G\le\e$ we have $\sum_{i\in F}\|x_i\|_{G}\le CM$. To show that the seminorm $\|\cdot\|$ on the group $G/H$ is  $C$-absolute, take any finite set $F$ and sequence $(y_i)_{i\in F}\in (G/H)^F$ such that $M:=\max_{A\subseteq F}\big\|\sum_{i\in F}y_i\|\le\e$. 
By the definition of the seminorm $\|\cdot\|$ on $G/H$, for every $i\in F$ there exists an element $x_i\in y_i$ such that $\|y_i\|\le\|x_i\|_G<\|y_i\|+\tfrac12\e$.

\begin{claim}\label{cl:quot} For every $k\in\{0,\dots,|F|\}$ and any set $E\subseteq F$ of cardinality $|E|=k$ we have
$$\big\|\sum_{i\in E}x_i\|_G=\big\|\sum_{i\in E}y_i\|.$$
\end{claim}

\begin{proof} The claim will be proved by induction on $k$. For $k=0$ and the empty set $E=\emptyset$ we have $$\big\|\sum_{i\in E}x_i\|_G=\|0_G\|_{G}=0=\|0_{G/H}\|=\big\|\sum_{i\in E}y_i\big\|.$$
Now assume that the equality in the claim has been proved for some nonnegative integer $k<|F|$. Take any set $E\subseteq F$ of cardinality $|E|=k+1$. Choose any index $j\in E$ and consider the set $E'=E\setminus\{j\}$. By the inductive assumption,
$  \big\|\sum_{i\in E'}x_i\big\|_G=\big\|\sum_{i\in E'}y_i\big\|\le M\le\e$. Assuming that $\big\|\sum_{i\in E}x_i\big\|_G\ne \big\|\sum_{i\in E}y_i\big\|$ and taking into account that $\sum_{i\in E}x_i\in\sum_{i\in E}y_i$, we conclude that $\big\|\sum_{i\in A}y_i\big\|<\|\sum_{i\in E}x_i\big\|_G$. Then we can find an element $h\in H$ such that $\big\|h+\sum_{i\in E}x_i\big\|_G<\big\|\sum_{i\in E}x_i\big\|_G$. The triangle inequality implies that
$$
\begin{aligned}
\|h\|_G&=\big\|h+\sum_{i\in E}x_i-\sum_{i\in E}x_i\big\|_G\le\big\|h+\sum_{i\in E}x_i\big\|_G+\big\|\sum_{i\in E}x_i\big\|_G<\\
&<\big\|\sum_{i\in E}x_i\big\|_G+\big\|\sum_{i\in E}x_i\|_G\le 2\big(\|x_j\|_G+\big\|\sum_{i\in E'}x_i\|_G\big)<\\
&<2\big(\|y_j\|+\tfrac12\e+\big\|\sum_{i\in E'}y_i\big\|\big)\le 2(M+\tfrac12\e+M)\le 5\e.
\end{aligned}
$$
Since $\{x\in H:\|x\|_G<5\e\}=\{0_H\}$, the element $h$ equals $0_H$ and we obtain a required contradiction:
$$\big\|\sum_{i\in E}x_i\big\|_G=\big\|h+\sum_{i\in E}x_i\big\|_G<
\big\|\sum_{i\in E}x_i\big\|_G,$$
showing that  $\big\|\sum_{i\in E}x_i\big\|_G= \big\|\sum_{i\in E}y_i\big\|$ and completing the proof of the claim.
\end{proof}

By Claim~\ref{cl:quot}, $$\max_{E\subseteq F}\big\|\sum_{i\in E}x_i\big\|_G=\max_{E\subseteq F}\big\|\sum_{i\in E}y_i\big\|=M\le \e.$$ Now the $C$-absoluteness of the seminorm $\|\cdot\|_G$ and the choice of $\e$ ensure that
$$\sum_{i\in F}\|y_i\|\le\sum_{i\in F}\|x_i\|_G\le CM,$$witnessing that the seminorm $\|\cdot\|$ on $G/H$ is $C$-absolute.
\end{proof}

Example~\ref{ex:R} and Proposition~\ref{p:quot} imply

\begin{lemma}\label{l:T2a} The seminorm $\|\cdot\|:\IR/\IZ\to[0,\infty)$, $\|\cdot\|:y\mapsto\min\{|x|:x\in y\}$, on the compact group $\IT=\IR/\IZ$  is $2$-absolute.
\end{lemma}

The following two lemmas can be easily deduced from the definitions.

\begin{lemma}\label{l:1a} Let $H$ be an open subgroup of a group $G$. If $\|\cdot\|_H:H\to[0,\infty)$ is a $C$-absolute seminorm on the group $H$, then the seminorm $$\|\cdot\|_G:G\to[0,\infty),\quad\|x\|_G=\begin{cases}\min\{1,\|x\|_H\}&\mbox{if $x\in H$};\\
1&\mbox{if $x\in G\setminus H$};
\end{cases}
$$on the group $G$ is $C$-absolute, too.
\end{lemma}

\begin{lemma}\label{l:ha} Let $h:H\to G$ be a continuous homomorphism between groups. If a seminorm $\|\cdot\|_G$ on the group $G$ is $C$-absolute for some positive real number $C$, then the seminorm
$$\|\cdot\|_H:H\to[0,\infty),\quad\|\cdot\|_H:x\mapsto\|h(x)\|_G,$$
is $C$-absolute.
\end{lemma}

\subsection{Locally Euclidean and locally compact groups}

A group $G$ is called {\em locally Euclidean} if some open neighborhood of $0_G$ in $G$ is homeomorphic to the Euclidean space $\IR^n$ for some $n\in\w$. In particular, every discrete group is locally Euclidean.

\begin{proposition}\label{p:lE} The topology of every locally Euclidean group $G$ is generated by a locally $2$-homogeneous absolute seminorm.
\end{proposition}

\begin{proof} By \cite[Example 75]{Pon}, the locally Euclidean group $G$ contains an open subgroup $H$, which is topologically isomorphic to the product $\IR^n\times\IT^k$ for some $n,k\in\w$. By Examples~\ref{ex:R2} and ~\ref{ex:R}, the norm $\|\cdot\|_\IR:\IR\to[0,\infty)$, $\|\cdot\|_\IR:x\mapsto |x|$, is locally $2$-homogeneous and absolute. By Corollary~\ref{cor:T2hom} and Lemma~\ref{l:T2a}, the seminorm $\|\cdot\|_\IT:\IT\to[0,\infty)$, $\|\cdot\|_\IT:y\mapsto \min\{|x|:x\in y\}$, is locally $2$-homogeneous and absolute. By Lemma~\ref{lem:prod2hom} and Corollary~\ref{c:prod-a}, the seminorms
$$\|\cdot\|_{\IR^n}:\IR^n\to[0,\infty),\quad\|\cdot\|_{\IR^n}:(x_i)_{i\in n}\mapsto\max_{i\in n}\|x_i\|_\IR,$$
$$\|\cdot\|_{\IT^k}:\IT^k\to[0,\infty),\quad\|\cdot\|_{\IT^k}:(x_i)_{i\in k}\mapsto\max_{i\in k}\|x_i\|_\IT,$$ and
$$\|\cdot\|:\IR^n\times\IT^k\to[0,\infty),\quad\|\cdot\|:(x,y)\mapsto\max\{\|x\|_{\IR^n},\|y\|_{\IT^k}\},$$
are locally $2$-homogeneous and absolute. It is clear that the topology of the group $\IR^n\times\IT^k$ is generated by the seminorm $\|\cdot\|$. Since $H$ is topologically isomorphic to $\IR^n\times\IT^k$, the topology of the group $H$ is generated by a locally $2$-homogeneous absolute seminorm $\|\cdot\|_H$. By Lemmas~\ref{l:1f} and \ref{l:1a}, the seminorm
$$\|\cdot\|_G:G\to[0,\infty),\quad\|\cdot\|_G:x\mapsto\begin{cases}\min\{1,\|x\|_H\}&\mbox{if $x\in H$};\\
1 &\mbox{if $x\in G\setminus H$};\\
\end{cases}
$$is locally $2$-homogeneous and absolute. It is clear that the seminorm $\|\cdot\|_G$ generates the topology of the group $G$.
\end{proof}

An actogroup $A\curvearrowright G$ is {\em locally compact} if the group $G$ is locally compact, that is every element of $G$ has a compact neighborhood in $G$.

\begin{proposition}\label{p:lc=>cas} The topology of any locally compact actogroup $A\curvearrowright G$ is generated by a family of contracting absolute seminorms.
\end{proposition}

\begin{proof} Let $U\subseteq G$ be a neighborhood of $0_G$ in $G$. Since $G$ is locally compact, we can replace $U$ by a smaller neighborhood, if necessary, and assume that $U$ has compact closure $\overline{U}$ in $G$.
By~\cite[Corollary 7.54]{HM}, there exists a surjective continuous homomorphism $h:G\to L$ onto a locally Euclidean group $L$ such that $h^{-1}[V]\subseteq U$ for some neighborhood $V\subseteq L$ of $0_L$. By Proposition~\ref{p:lE}, the topology of the locally Euclidean group $L$ is generated by a locally $2$-homogeneous absolute seminorm $\|\cdot\|_L:L\to[0,\infty)$. So, $\{x\in L:\|x\|_L\le\e\}\subseteq V$ for some $\e>0$. By Lemmas~\ref{lem:hom2hom} and \ref{l:ha}, the seminorm $$\|\cdot\|_G:G\to[0,\infty),\quad \|\cdot\|_G:x\mapsto \|h(x)\|_L$$ on the  group $G$ is locally $2$-homogeneous and absolute. Since the closed neighborhood $$\{x\in G:\|x\|_G\le\e\}\subseteq \{x\in G:\|h(x)\|_L\le\e\}\subseteq\{x\in G:h(x)\in V\}\subseteq U\subseteq\overline U$$is compact, we can apply Theorem~\ref{t:Ravsky} and conclude that the seminorm $\|\cdot\|_G$ is contracting.
\end{proof}

\subsection{Obsolete seminorms}

\begin{definition} A seminorm $\|\cdot\|:G\to[0,\infty)$ on a group $G$ is defined to be
{\em obsolete} if for every unconditionally convergent series $\sum_{i\in\w}x_i$ in $G$ the series $\sum_{i\in\w}\|x_i\|$ converges.
\end{definition}

\begin{lemma}\label{l:a=>o} Each absolute seminorm $\|\cdot\|$ on a group $G$ is obsolete.
\end{lemma}

\begin{proof} Being absolute, the seminorm $\|\cdot\|$ is $C$-absolute for some positive real number $C$. This means that there exists $\e>0$ such that such that for any finite set $F$ and sequence $(x_i)_{i\in F}\in G^F$ with $M:=\max_{E\subseteq F}\|\sum_{i\in E}x_i\|\le \e$ we have $\sum_{i\in F}\|x_i\|\le CM$.

 To show that $\|\cdot\|$ is obsolete, fix any unconditionally convergent series $\sum_{i\in\w}x_i$ in $G$. The unconditional convergence of $\sum_{i\in\w}x_i$ yields a finite set $K\subseteq\w$ such that $\big\|\sum_{i\in F}x_i\big\|<\e$ for any finite set $F\subset\w\setminus K$.

 Suppose for a contradiction that the series $\sum_{i\in\w}\|x_i\|$ diverges in the real line. Then we can find a finite set $F\subset\w\setminus K$ such that  $\sum_{i\in F}\|x_i\|>C\e$.
 The choice of the sets $K$ and $F$ ensures that $$M:=\max_{E\subseteq F}\big\|\sum_{i\in E}x_i\big\|<\e$$ and hence $\sum_{i\in F}\|x_i\|\le CM\le C\e$, a contradiction.
 \end{proof}

The obsoleteness is preserved by standard operations over seminorms.

\begin{lemma}\label{l:ho} Let $h:H\to G$ be a continuous homomorphism between groups. If $\|\cdot\|_G$ is an obsolete seminorm on the group $G$, then the seminorm
$$\|\cdot\|_H:H\to[0,\infty),\quad \|\cdot\|_H:x\mapsto\|h(x)\|_G,$$
is an obsolete seminorm on the group $H$.
\end{lemma}

\begin{proof} To prove that the seminorm $\|\cdot\|_H$ is obsolete, take any unconditionally convergent series $\sum_{i\in\w}x_i$ in the group $H$. The continuity of the homomorphism $h$ implies that $\sum_{i\in\w}h(x_i)$ is an unconditionally convergent series in the group $G$. Since the seminorm $\|\cdot\|_G$ is obsolete, the series $\sum_{i\in\w}\|h(x_i)\|_G=\sum_{i\in\w}\|x_i\|_H$ converges, witnessing that the seminorm $\|\cdot\|_H$ is obsolete.
\end{proof}

\begin{lemma}\label{l:prod-o} For two obsolete seminorms $\|\cdot\|_G:G\to[0,\infty)$ and $\|\cdot\|_H:H\to[0,\infty)$ on groups $G,H$, the seminorm
$$\|\cdot\|:G\times H\to[0,\infty),\quad \|\cdot\|:(g,h)\mapsto\max\{\|g\|_G,\|h\|_H\},$$
on $G\times H$ is obsolete.
\end{lemma}

\begin{proof} Let $\sum_{i\in\w}z_i$ be any unconditionally convergent series in the group $G\times H$. For every $i\in\w$ write the element $z_i\in G\times H$ as a pair $(x_i,y_i)$ for some $x_i\in G$ and $y_i\in H$. The unconditional convergence of the series $\sum_{i\in\w}z_i$ implies the unconditional convergence of the series $\sum_{i\in\w}x_i$ and $\sum_{i\in\w}y_i$. By the obsoleteness of the seminorms $\|\cdot\|_G$ and $\|\cdot\|_H$, the series $\sum_{i\in\w}\|x_i\|_G$ and $\sum_{i\in\w}\|y_i\|_H$ converge in the real line. Then
$$\sum_{i\in\w}\|z_i\|=\sum_{i\in\w}\max\{\|x_i\|_G,\|y_i\|_H\}\le\sum_{i\in\w}(\|x_i\|_G+\|y_i\|_H)\le\sum_{i\in\w}\|x_i\|_G+\sum_{i\in\w}\|y_i\|_H<\infty,$$
which means that the series $\sum_{i\in\w}\|z_i\|$ converges, and therefore the seminorm $\|\cdot\|$ on $G\times H$ is obsolete.
\end{proof}

\subsection{Unconditional seminorms on actogroups}

\begin{definition}  A seminorm $\|\cdot\|:G\to[0,\infty)$ on an actogroup $A\curvearrowright G$  is called
\begin{itemize}
\item {\em $C$-unconditional} for some real number $C$ if there exist a positive real number $\e$ and a neighborhood $V\subseteq A$ of $\mathsf 0_A$ such that for any finite set $F$ and sequences $(x_i)_{i\in F}\in G^F$ and $(a_i)_{i\in F}\in V^F$ with $M:=\max\limits_{E\subseteq F}\big\|\!\sum\limits_{i\in E}x_i\big\|\le \e$ we have $\big\|\sum\limits_{i\in F}a_ix_i\big\|\le CM$;
\item {\em unconditional} if $\|\cdot\|$ is $C$-unconditional for some positive real number $C$.
\end{itemize}
\end{definition}

\begin{lemma}\label{l:aL=>u}  If a seminorm $\|\cdot\|$ on an actogroup $A\curvearrowright G$ is absolute and Lipschitz, then $\|\cdot\|$ is unconditional.
\end{lemma}

\begin{proof} Assuming that $\|\cdot\|$ is absolute and Lipschitz, find positive real numbers $L,C,\e$ and a neighborhood $V\subseteq A$ of $\mathsf 0_A$ such that
\begin{itemize}
\item $\sum_{i\in F}\|x_i\|\le C\cdot\max_{E\subseteq F}\big\|\sum_{i\in E}x_i\big\|$ for any finite set $F$ and sequence $(x_i)_{i\in F}\in G^F$ with $\max_{E\subseteq F}\big\|\sum_{i\in E}x_i\big\|\le\e$;
\item $\|ax\|\le L\|x\|$ for any $a\in V$ and $x\in G$ with $\|x\|\le\e$.
\end{itemize}
Then for any finite set $F$ and sequences $(x_i)_{i\in F}\in G^F$ and $(t_i)_{i\in F}\in V^F$ with $\max\limits_{E\subseteq F}\big\|\sum\limits_{i\in E}x_i\big\|\le \e$ we have
$$\big\|\sum_{i\in F}t_ix_i\big\|\le\sum_{i\in F}\|t_ix_i\|\le \sum_{i\in F}L\|x_i\|\le LC\max_{E\subseteq F}\big\|\sum_{i\in E}x_i\big\|,$$
which means that the seminorm $\|\cdot\|$ is $LC$-unconditional and hence unconditional.
\end{proof}

\subsection{Absolutely and unconditionally seminormable actogroups}

\begin{definition}\label{d:aou} An actogroup $A\curvearrowright G$ is called
\begin{enumerate}
\item  {\em absolutely seminormable}   if the topology of $G$ is generated by a family of Lipschitz  absolute seminorms;
\item  {\em obsoletely seminormable}   if the topology of $G$ is generated by a family of Lipschitz  obsolete seminorms;
\item  {\em unconditionally seminormable} if the topology of $G$ is generated by a family of unconditional seminorms.
\end{enumerate}
\end{definition}

Lemmas~\ref{l:a=>o} and \ref{l:aL=>u} imply

\begin{lemma}\label{l:as=>ous} Every absolutely seminormable actogroup is obsoletely and unconditionally seminormable.
\end{lemma}

\begin{theorem}\label{t:uo-H} A actogroup $A\curvearrowright G$ is Hirsch if one of the following conditions holds:
\begin{enumerate}
\item $A\curvearrowright G$ is unconditionally seminormable;
\item $A\curvearrowright G$ is obsoletely seminormable;
\item $A\curvearrowright G$ is absolutely seminormable.
\end{enumerate}
\end{theorem}

\begin{proof}   To prove that the actogroup $A\curvearrowright G$ is Hirsch, fix any unconditionally convergent series $\sum_{n\in\w}x_i$ in $G$ and any neighborhood $U$ of $0_G$ in $G$.
\smallskip

1. First we assume that $A\curvearrowright G$ is unconditionally seminormable. Then there exists an unconditional seminorm $\|\cdot\|$ on $G$ such that $\{x\in G:\|x\|\le 2 \e\}\subseteq U$ for some $\e>0$. Since $\|\cdot\|$ is unconditional, there exist positive numbers $C,\delta$ and a neighborhood $V\subseteq A$ of $\mathsf 0_A$ such that for every finite set $F$ and any sequences $(y_i)_{i\in F}\in G^F$ and $(a_i)_{i\in F}\in V^F$ with $M=\max_{E\subseteq F}\big\|\sum_{i\in E}y_i\|\le\delta$, we have $\big\|\sum_{i\in F}a_iy_i\|\le CM$.

Since the series $\sum_{i\in\w}x_i$ is unconditionally convergent, there exists a finite set $E\subseteq\w$ such that $\big\|\sum_{i\in F}x_i\big\|<\min\{\delta,\tfrac\e{C}\}$ for every finite set $F\subseteq\w\setminus E$.

By the continuity of the action $A\times G\to G$, for every $i\in E$, there exists a neighborhood $V_i\subseteq A$ of $\mathsf 0_A$ such that $\sup_{a\in V_i}\|ax_i\|\le \frac{\e}{|E|}$. Replacing the neighborhood $V$ by the neighborhood $V\cap\bigcap_{i\in E}V_i$, we can additionally assume that $\max_{i\in E}\sup_{a\in V}\|ax_i\|\le\frac\e{|E|}$.

Now observe that for every finite set $F$ and sequence $(a_i)_{i\in F}\in V^F$ we have
\begin{multline*}
\big\|\sum_{i\in F}a_ix_i\big\|\le\sum_{i\in F\cap E}\|a_ix_i\|+\big\|\sum_{i\in F\setminus E}a_ix_i\big\|\le\sum_{i\in F\cap E}\frac{\e}{|E|}+C\max_{K\subseteq F\setminus E}\big\|\sum_{i\in K}x_i\big\|\le\\
 \le \frac{\e|F\cap E|}{|E|}+C\frac\e{C}\le2\e
\end{multline*}
and hence $\sum_{i\in F}a_ix_i\in U$.
\smallskip

2. Next, assume that $A\curvearrowright G$ is obsoletely seminormable. Then there exists a Lipschitz obsolete seminorm $\|\cdot\|$ on $G$ such that $\{x\in G:\|x\|\le 2\e\}\subseteq U$ for some $\e>0$. By the obsoleteness of $\|\cdot\|$, we have $\sum_{n\in\w}\|x_i\|<\infty$. Since the seminorm $\|\cdot\|$ is Lipschitz, there exist positive numbers $L\ge 1$, $\e'\le\e$
 and a neigborhood $V'\subseteq A$ of $\mathsf 0_A$ such that $\|ax\|\le L\|x\|$ for every $a\in V'$ and $x\in G$ with $\|x\|\le\e'$.  Since $\sum_{i\in\w}\|x_i\|<\infty$, there exists a nonempty finite set $E\subset\w$ such that $\sum_{i\in\w\setminus E}\|x_i\|\le\frac{\e'}{L}$. For every $i\in \w\setminus E$ we have $\|x_i\|\le\sum_{j\in\w\setminus E}\|x_j\|\le\e'$
   and hence $\sup_{a\in V'}\|ax_i\|\le L\|x_i\|$.  By the continuity of the seminorm $\|\cdot\|$ and continuity of the action $A\times G\to G$, for every $i\in E$, there exists a neighborhood $V_i\subseteq A$ of $\mathsf 0_A$ such that $\sup_{a\in V_i}\|ax_i\|\le \frac{\e}{|E|}$. We claim that the neighborhood $V=V'\cap \bigcap_{i\in I}V_i$ witnesses that $G$ is Hirsch. Take any finite set $F\subseteq \w$ and any sequence $(a_i)_{i\in F}\in V^F$.
It follows that
\begin{multline*}
\big\|\sum_{i\in F}a_ix_i\big\|\le\sum_{i\in F\cap E}\|a_ix_i\|+\sum_{i\in F\setminus E}\|a_ix_i\|\le\sum_{i\in F\cap E}\frac{\e}{|E|}+\sum_{i\in F\setminus E}L\|x_i\|\le
 \frac{\e|F\cap E|}{|E|}+L\frac\e{L}\le 2\e
\end{multline*}
and hence $\sum_{i\in F}a_ix_i\in U$.
\smallskip

3. If $A\curvearrowright G$ is absolutely seminormable, then by Lemma~\ref{l:as=>ous}, $A\curvearrowright G$ is both obsoletely seminormable and unconditionally seminormable and hence $A\curvearrowright G$ is Hirsch by any of the preceding cases.
\end{proof}

Lemmas~\ref{l:a=>o}, \ref{l:aL=>u} and Theorem~\ref{t:uo-H} imply that for any actogroup we have the implications.
$$
\xymatrix{
\mbox{absolutely seminormable}\ar@{=>}[r]\ar@{=>}[d]&\mbox{obsoletely seminormable}\ar@{=>}[d]\\
\mbox{unconditionally seminormable}\ar@{=>}[r]&\mbox{Hirsch}
}
$$

Since every contracting seminorm of an actogroup is Lipschitz, Proposition~\ref{p:lc=>cas} and Theorem~\ref{t:uo-H} imply:

\begin{theorem}\label{t:lc-H} Every locally compact actogroup is absolutely seminormable and  Hirsch.
\end{theorem}





\subsection{Tychonoff products of actogroups}

\begin{definition} The {\em Tychonoff product} $\prod_{i\in I}(A_i\curvearrowright G_i)$ of a nonempty family of actogroups $(A_i\curvearrowright G_i)_{i\in I}$ is the group $\prod_{i\in I}G_i$ endowed with the action
$$\prod_{i\in I}A_i\times\prod_{i\in I}G_i\to\prod_{i\in I}G_i,\quad \big((a_i)_{i\in I},(x_i)_{i\in I}\big)\mapsto (a_ix_i)_{i\in I}$$of the topological space $\prod_{i\in I}A_i$ whose distinguished point is the unique point of the singleton $\prod_{i\in I}\{\mathsf 0_{A_i}\}$.

If all the actogroups $A_i\curvearrowright G_i$ are equal to some fixed actogroup $A\curvearrowright G$, then the product $\prod_{i\in I}(A_i\curvearrowright G_i)$ will be denoted by $(A\curvearrowright G)^I$.
\end{definition}

\begin{theorem} The Tychonoff product $\prod_{i\in I}(A_i\curvearrowright G_i)$ of a nonempty family of absolutely (resp. obsoletely) seminormable actogroups is absolutely (resp. obsoletely) seminormable and Hirsch.
\end{theorem}

\begin{proof} For a subset $F\subseteq I$ let $G_F=\prod_{i\in F}G_i$ and $\pr_F:G\to G_F$, $\pr_F:(x_i)_{i\in I}\mapsto (x_i)_{i\in F}$, be the projection homomorphism. It follows that $\prod_{i\in I}G_i=G_I$. Given any neighborhood $U\subseteq G_I$ of $0_{G_I}$, find a finite set $F\subseteq I$ and a neighborhood $W$ of $0_{G_F}$ in the group $G_F$ such that $\pr_F^{-1}[W]\subseteq U$. For every $i\in F$ find a neighborhood $W_i$ of $0_{G_i}$ in $G_i$ such that $\prod_{i\in F}W_i\subseteq W$. Since the actogroup $A_i\curvearrowright G_i$ is absolutely (resp. obsoletely) seminormable, there exists a Lipschitz absolute (resp. obsolete) seminorm $\|\cdot\|_i:G_i\to[0,\infty)$ such that $\{x\in G_i:\|x\|_i\le\e_i\}\subseteq W_i$ for some $\e_i>0$. By Lemma~\ref{l:prod-a} (resp. \ref{l:prod-o}), the seminorm
$$\|\cdot\|_F:G_F\to[0,\infty),\quad\|\cdot\|_F:(x_i)_{i\in F}\mapsto\max_{i\in F}\|x\|_i,$$is absolute (resp. obsolete). By Lemma~\ref{l:ha} (resp. \ref{l:ho}), the seminorm
$$\|\cdot\|:G_I\to [0,\infty),\quad \|\cdot\|:x\mapsto\|\pr_F(x)\|_F,$$
is absolute (resp. obsolete). It is clear that $$\{x\in G_I:\|x\|\le\min_{i\in F}\e_i\}\subseteq \pr_F^{-1}\Big[\prod_{i\in F}W_i\Big]\subseteq \pr_F^{-1}[W]\subseteq U.$$ It remains to prove that the seminorm $\|\cdot\|$ is Lipschitz. For every $i\in I$ the  Lipschitz property of the seminorm $\|\cdot\|_i$ yields positive numbers $L_i,\delta_i$ and a neighborhood $V_i\subseteq A_i$ of $\mathsf 0_{A_i}$ such that $\sup_{a\in V_i}\|ax\|_i\le L_i\|x\|_i$ for any $x\in G_i$ with $\|x\|_i\le\delta_i$. Let $L=\max_{i\in F}L_i$, $\delta=\min_{i\in F}\delta_i$ and $V=\bigcap_{i\in F}\{(a_j)_{j\in I}\in \prod_{j\in I}A_j:a_i\in V_i\}$. It is easy to see that for every $a=(a_j)_{j\in I}\in V$ and $x=(x_j)_{j\in I}\in G$ with $\|(x_j)_{j\in I}\|=\max_{i\in F}\|x_i\|_i\le\delta=\min_{i\in F}\delta_i$ we have
$$\|ax\|=\max_{i\in F}\|a_ix_i\|_i\le \max_{i\in F}L_i\|x_i\|_i\le L\max_{i\in F}\|x_i\|_i=L\|x\|,$$
which means that the seminorm $\|\cdot\|$ is Lipschitz. Therefore, the actogroup $\prod_{i\in I}(A_i\curvearrowright G_i)$ is absolutely (resp. obsoletely) seminormable.
By Theorem~\ref{t:uo-H}, this actogroup is Hirsch.
\end{proof}

\begin{corollary} For every nonempty set $I$ and every absolutely (resp. obsoletely) seminormable actogroup $A\curvearrowright G$ the actogroup $(A\curvearrowright G)^I$ is absolutely (resp. obsoletely) seminormable and Hirsch.
\end{corollary}

\begin{theorem} For any nonempty set $I$ and any locally compact actogroup $A\curvearrowright G$, the actogroup $(A\curvearrowright G)^I$ is absolutely seminormable and Hirsch.
\end{theorem}



\subsection{Function actogroups}

A family $\K$ of compact subsets of a topological space $X$ is called an {\em ideal of compact sets} if for any sets $A,B\in\mathcal K$ and any closed subset $K\subseteq A\cup B$ we have $K\in\mathcal K$.

For topological spaces $X,Y$ and an ideal $\K$ of compact subsets of $X$, let $C_\K(X,Y)$ be the space of continuous functions from $X$ to $Y$, endowed with the topology generated by the subbase consisting of the sets
$$[K;U]=\{f\in C_\K(X,Y):f(K)\subseteq U\}$$ where $K\in\K$ and $U$ is an open subset of $Y$.

If $\K$ is the ideal of all compact (finite) subsets of $X$, then the space $C_\K(X,Y)$ is denoted by $C_k(X,Y)$ (resp. $C_p(X,Y)$).

If $X$ is a discrete space then  $C_p(X,Y)$ is equal to the power $Y^X$, endowed with the Tychonoff product topology.

Let $X$ be a topological space and $\K$ be an ideal of compact Hausdorff subspaces of $X$. Given topological spaces $Y,Z,T$ and a continuous map $\alpha:Y\times Z\to T$, consider the map  $\tilde \alpha:C_\K(X,Y)\times C_\K(X,Z)\to C_\K(X,T)$ assigning to every pair of functions $(f,g)\in C_\K(X,Y)\times C_\K(X,Z)$ the function $\tilde \alpha(f,g)\in C_\K(X,T)$ defined by $\tilde \alpha(f,g)(x)=\alpha(f(x),g(x))$ for $x\in X$. Since the family $\K$ consists of compact Hausdorff subspaces of $X$, we can apply Proposition 4.2 of \cite{BL} and conclude that the map $\tilde\alpha$ is continuous.

This implies that for every actogroup $A\curvearrowright G$ we can consider the actogroup $C_\K(X,A)\curvearrowright C_\K(X,G)$ endowed with the action $$C_\K(X,A)\times C_\K(X,G)\to C_\K(X,G),\quad (a,g)\mapsto ag,$$ where $ag:x\mapsto a(x)g(x)$. The distinguished point of the topological space $C_\K(X,A)$ is the constant function $X\to\{\mathsf 0_A\}\subseteq A$.

\begin{theorem}\label{t:u-CK} Let $X$ be a topological space and $\K$ be an ideal of compact Hausdorff
 subspaces of $X$. If $A\curvearrowright G$ is an unconditionally seminormable actogroup, then the actogroup $C_\K(X,A)\curvearrowright C_\K(X,G)$ is unconditionally seminormable and  Hirsch.
\end{theorem}

\begin{proof} Given any neighborhood $\U$ of zero in the group $C_\K(X,G)$, find a set $K\in\K$ and a neighborhood $U$ of zero in $G$ such that $[K,U]\subseteq\U$. Since the actogroup $A\curvearrowright G$ is unconditionally seminormable, there exists an unconditional seminorm $\|\cdot\|_G$ on $G$ such that $\{x\in G:\|x\|_G\le\e\}\subseteq U$ for some $\e>0$. By the unconditionality of $\|\cdot\|_G$, there exist positive real numbers $C,\delta$ and a neighborhood $V\subseteq A$ of zero $\mathsf 0_A$ such that for any finite set $F$ and any sequences $(x_i)_{i\in F}\in G^F$ and $(a_i)_{i\in F}\in V^F$ with $\max_{E\subseteq F}\big\|\sum_{i\in E}x_i\big\|_G\le\delta$ we have $\big\|\sum_{i\in F}a_ix_i\big\|_G\le C\cdot\max_{E\subseteq F}\big\|\sum_{i\in E}x_i\big\|_G$.

Consider the seminorm $$\|\cdot\|_K:C_\K(X,G),\;\|\cdot\|_K:f\mapsto \|f\|_K=\sup_{x\in K}\|f(x)\|_G,$$
on the group $C_\K(X,G)$.
The compactness of $K$ implies that the seminorm $\|\cdot\|_K$ is well-defined and continuous. Observe that
$$\{f\in C_\K(X,G):\|f\|_K\le\e\}\subseteq [K;U]\subseteq\U.$$
It remains to show that the seminorm $\|\cdot\|_K$ is unconditional. The definition of the topology on the space $C_\K(X,T)$ ensures that $[K;V]=\{f\in C(X,A):f[K]\subseteq V\}$ is a neighborhood of zero in the acting topological space $C_\K(X,A)$.

Now take any finite set $F$ and sequences $(f_i)_{i\in F}\in C_\K(X,G)^F$ and $(a_i)_{i\in F}\in [K;V]^F$ such that $M:=\max_{E\subseteq F}\big\|\sum_{i\in E}f_i\|_K\le\e$. Observe that for every $x\in K$ we have $\{a_i(x)\}_{i\in F}\subseteq V$ and $$\max_{E\subseteq F}\big\|\sum_{i\in E}f_i(x)\big\|_G\le \max_{E\subseteq F}\big\|\sum_{i\in E}f_i\big\|_K=M\le\e$$and hence
$$\big\|\sum_{i\in F}a_i(x)f_i(x)\big\|_G\le C\cdot\max_{E\subseteq F}\big\|\sum_{i\in E}f_i(x)\big\|_{G}\le CM.$$
Then $$\big\|\sum_{i\in F}a_if_i\big\|_K=\max_{x\in K}\big\|\sum_{i\in F}a_i(x)f_i(x)\big\|_{G}\le CM=C\max_{E\subseteq F}\big\|\sum_{i\in E}f_i\|_K,$$
witnessing that the seminorm $\|\cdot\|_K$ on $C_\K(X,G)$ is unconditional.

Therefore, the topology of the group $C_\K(X,G)$ is generated by unconditional seminorms and by Theorem~\ref{t:uo-H}, the actogroup $C_\K(X,A)\curvearrowright C_\K(X,G)$ is Hirsch.
\end{proof}

\begin{theorem}\label{t:lc-CK} Let $X$ be a topological space and $\K$ be an ideal of compact Hausdorff subspaces of $X$. For every locally compact actogroup $A\curvearrowright G$, the actogroup $C_\K(X,A)\curvearrowright C_\K(X,G)$ is unconditionally seminormable and  Hirsch.
\end{theorem}

\begin{proof} By Theorem~\ref{t:lc-H} and Lemma~\ref{l:as=>ous}, the actogroup $A\curvearrowright G$ is absolutely seminormable and unconditionally seminormable.  By Theorem~\ref{t:u-CK}, the actogroup $C_\K(X,A)\curvearrowright C_\K(X,G)$ is unconditionally seminormable and  Hirsch.
\end{proof}

\section{Applications to Hirsch topological rings}

In this section we apply the results of the preceding section to topological rings. First we write down Definition~\ref{d:aou}(3) for the partial case of topological rings.

\begin{definition} A topological ring $R$ is {\em unconditionally seminormable} if for every neighborhood $U\subseteq R$ of $0_R$ there exist a seminorm $\|\cdot\|:R\to[0,\infty)$ on the additive group of $R$, a neighborhood $V\subseteq R$ of $0_R$, and positive real numbers $C,\e$ such that $\{x\in R:\|x\|<\e\}\subseteq U$ and for any finite set $F$ and sequences $(x_i)_{i\in F}\in R^F$ and $(a_i)_{i\in F}\in V^F$ with $M:=\max_{E\subseteq F}\big\|\sum_{i\in E}x_i\big\|\le\e$ we have $$\big\|\sum_{i\in F}a_ix_i\big\|\le CM.$$
\end{definition}

The following theorems are partial cases of Theorems~\ref{t:uo-H}, \ref{t:lc-H}, \ref{t:u-CK} and \ref{t:lc-CK}, respectively.

\begin{theorem}\label{t:5.2} Every unconditionally seminormable topological ring is Hirsch.
\end{theorem}

\begin{theorem}\label{t:lc-Hr} Every locally compact topological ring is unconditionally seminormable and Hirsch.
\end{theorem}


\begin{theorem}\label{t:CK-ring}  Let $X$ be a topological space, $\K$ be an ideal of compact Hausdorff subspaces of $X$, and $R$ be a topological ring which is either unconditionally seminormable or locally compact. Then the topological ring $C_\K(X,R)$ is unconditionally seminormable and Hirsch.
\end{theorem}

Theorem~\ref{t:CK-ring} has two corollaries.

\begin{corollary} For every locally compact topological ring $R$ and any nonempty set $I$ the topological ring $R^I$ is unconditionally seminormable and Hirsch.
\end{corollary}

\begin{corollary}\label{c:5.7} For every Hausdorff topological space $X$ and any locally compact topological ring $R$ the topological rings $C_p(X,R)$ and $C_k(X,R)$ are unconditionally seminormable and Hirsch.
\end{corollary}

\section{Banach actospaces}

A {\em Banach actospace} $A\curvearrowright X$ is a Banach space $X$ endowed with a continuous bilinear action $A\times X\to X$, $(a,x)\mapsto ax$, of a Banach space $A$.  All Banach spaces considered in this paper are over the field $\mathbb F$ of real or complex numbers. It is clear that each Banach actospace is an actogroup. Banach actospaces with the Hirsch property admit a nice characterization.

\begin{theorem}\label{t:Banact} For a  Banach actospace $A\curvearrowright X$ the following conditions are equivalent:
\begin{enumerate}
\item $A\curvearrowright X$ is Hirsch;
\item the norm of $X$ is unconditional;
\item there exists a positive real constant $C$ such that for every finite set $F$ and sequences $(a_n)_{n\in F}\in A^F$ and $(x_n)_{n\in F}\in X^F$ we have
$\|\sum\limits_{n\in F}a_nx_n\|\le C\cdot\max\limits_{n\in F}\|a_n\|\cdot\max\limits_{E\subseteq F}\|\sum\limits_{n\in E}x_n\|$;
\item for any unconditionally convergent series $\sum_{n\in\w}x_n$ in $X$ and any bounded sequence $(a_n)_{n\in\w}$ in $A$, the series $\sum_{n\in\w}a_nx_n$ converges (unconditionally) in $X$.
\end{enumerate}
\end{theorem}

\begin{proof} We shall prove the implications $(3)\Ra(2)\Ra(1)\Ra(3)\Leftrightarrow(4)$.
\smallskip

$(3)\Ra(2)$ Assume that the condition (3) is satisfied for some positive real constant $C$. We claim that the norm of the Banach space $X$ is $C$-uncounditional. Given a positive real number $\e$, consider the neighborhood $V=\{a\in A:\|a\|<1\}$ of zero in the Banach space $A$. By the choice of $C$, for any finite set $F$ and sequences $(a_n)_{n\in F}\in V^F$ and $(x_n)_{n\in F}\in X^F$ we have
$$\big\|\sum_{n\in F}a_nx_n\big\|\le C\max_{n\in F}\|a_i\|\cdot\max_{E\subseteq F}\big\|\sum_{n\in E}x_n\big\|\le C\max_{E\subseteq F}\big\|\sum_{n\in E}x_n\big\|,$$
witnessing that the norm of $X$ is $C$-unconditional.
\smallskip

The implication $(2)\Ra(1)$ follows from Theorem~\ref{t:uo-H}(1).
\smallskip

$(1)\Ra(3)$ Assume that the Banach actogroup $A\curvearrowright X$ is Hirsch. Consider the linear space $\Sigma X$ of all sequences $(x_n)_{n\in\w}$ in $X$ that have finite norm
$$\|(x_n)_{n\in\w}\|_{\Sigma X}=\sup_{F\in[\w]^{<\w}}\big\|\sum_{n\in F}x_n\big\|,$$
where $[\w]^{<\w}$ is the family of finite subsets of $\w$.
It is easy to see that $\Sigma X$ is a Banach space with respect to the norm $\|\cdot\|_{\Sigma X}$. Let $\Sigma_0 X$ be the (closed) subspace of $\Sigma X$ consisting of sequences $(x_n)_{n\in\w}$ for which the series $\sum_{n\in\w}x_n$ converges unconditionally in $X$. Let $\ell_\infty[A]$ be the Banach space of all bounded sequences $(a_n)_{n\in\w}$ in $A$, endowed with the norm $\|(a_n)_{n\in\w}\|_{\ell_\infty[A]}=\sup_{n\in\w}\|a_n\|_A$. By the Hirsch property of $X$, the bilinear operator
$$T:\ell_\infty[A]\times \Sigma_0X\to\Sigma X,\quad T:\big((a_n)_{n\in\w},(x_n)_{n\in\w}\big)\mapsto (a_nx_n)_{n\in\w},$$is well-defined. The operator $T$ has closed graph
$$\big\{\big((a_n)_{n\in\w},(x_n)_{n\in\w},(y_n)_{n\in\w}\big)\in \ell_\infty[A]\times\Sigma_0X\times\Sigma X:\forall n\in\w\;\;(a_nx_n=y_n)\big\}$$and hence is continuous by the Closed Graph Theorem of Fernandez \cite{Fern}. Then $T$ has bounded norm
$$\|T\|=\sup\{\|T(a,x)\|_{\Sigma X}:(a,x)\in\ell_\infty[A]\times\Sigma_0 X,\;\;\max\{\|a\|_{\ell_\infty[A]},\|x\|_{\Sigma_0 X}\}\le 1\}$$and hence
$$\big\|\sum_{n\in\w}a_nx_n\big\|\le \|T\|\cdot\sup_{n\in\w}\|a_n\|\cdot \sup_{F\in[\w]^{<\w}}\big\|\sum_{n\in F}x_n\big\|,$$for any sequences $(a_n)_{n\in\w}\in\ell_\infty[A]$ and $(x_n)_{n\in\w}\in\Sigma_0 X$. The latter inequality implies that the condition (3) holds for the constant $C=\|T\|$.
\smallskip

The equivalence $(3)\Leftrightarrow(4)$ is proved in Theorem 3.2 of \cite{BanKad}.
\end{proof}

The following theorem follows from Theorem~\ref{t:Banact} and Theorem 4.1 in \cite{BanKad}.

\begin{theorem} A Banach actospace $A\curvearrowright X$ is Hirsch if $A$ is a Hilbert space possessing an orthonormal basis $(e_n)_{n\in\w}$ such that for every $x\in X$ the series $\sum_{n\in\w}e_nx$ converges unconditionally in $X$.
\end{theorem}

For a number $p\in[1,\infty]$ let $\ell_p$ be the Banach space of  all functions $x:\w\to\IF$ to the field $\IF$ of real or complex numbers such that the norm
$$\|x\|_{\ell_p}=\begin{cases}\big(\sum_{n\in\w}|x(n)|^p\big)^{\frac1p}&\mbox{if $p<\infty$;}\\
\phantom{i}\sup_{n\in\w}|x(n)|&\mbox{if $p=\infty$;}
\end{cases}
$$is finite. For any $p,q\in[1,\infty]$ and functions $x\in\ell_p$ and $y\in\ell_q$, the function $xy:\w\to\IF$, $xy:n\mapsto x(n)y(n)$, belongs to $\ell_q$, witnessing that the Banach actospace $\ell_p\curvearrowright\ell_q$ is well-defined.
The Hirsch properties of this actospace are described in the following theorem that follows from Theorem~\ref{t:Banact} and Theorem 1.4 in \cite{BanKad}.

\begin{theorem}\label{t:lpq} For numbers $p,q\in[1,\infty]$ the Banach actospace  $\ell_p\curvearrowright \ell_q$ is
\begin{enumerate}
\item Hirsch if $p\in[1,2]$ or $q=\infty$;
\item not Hirsch if $p>\max\{2,q\}$.
\end{enumerate}
\end{theorem}

\section{Banach rings}

In this section we recognize Hirsch rings among Banach rings and also present examples of Banach rings which are not Hirsch.

A {\em Banach ring} is a Banach space $X$ endowed with a binary operation $\cdot:X\times X\to X$ turning $X$ into a topological ring.

\subsection{Hirsch Banach rings}
For a Banach ring $R$ and a compact Hausdorff space $K$, let $C(K,R)$ be the Banach ring endowed with the norm $\|\cdot\|_K:C(K,R)\to[0,\infty)$, $\|\cdot\|_K:f\mapsto \sup_{x\in K}\|f(x)\|$, where $\|\cdot\|$ is the norm of the Banach space $R$.

\begin{theorem}\label{t:B-CK} For every finite-dimensional Banach ring $R$ and compact Hausdorff space $K$, the Banach ring $C(K,R)$ is Hirsch.
\end{theorem}

\begin{proof} Being finite-dimensional, the Banach space $R$ is locally compact. Now Corollary~\ref{c:5.7} implies that the Banach ring $C(K,R)=C_k(K,R)$ is Hirsch.
\end{proof}

By $\IR,\IC,\IH$ we denote the Banach rings of real, complex and quaternion numbers, respectively. Theorem~\ref{t:B-CK} implies

\begin{corollary}\label{c:6.2} For any compact Hausdorff space $K$ the Banach rings $C(K,\IR)$, $C(K,\IC)$ and $C(K,\IH)$ are Hirsch.
\end{corollary}

Let $\ell_\infty$ be the Banach ring of all bounded sequences of real (or complex) numbers,
endowed with the sup-norm $\|(x_i)_{i\in\w}\|_{\ell_\infty}=\sup_{i\in\w}|x_i|$ and the operation of coordinatewise multiplication.  Let $c$ be the Banach subring of $\ell_\infty$ consisting of convergent sequences and  $c_0$ be the Banach subring of $c$, consisting of sequences that converge to zero.

\begin{theorem}\label{t:6.3} The Banach rings $c_0$, $c$, and $\ell_\infty$ are Hirsch.
\end{theorem}

\begin{proof} Let $\beta\w$ be the Stone-\v Cech compactification of the discrete space $\w$ of finite ordinals and $\IF$ be the topological field of real (or complex) numbers. It is well-known (and easy to see) that the restriction operator $C(\beta\w,\IF)\to\ell_\infty$, $f\mapsto f{\restriction}_\w$, is a topological isomorphism of the Banach rings $C(K,\IF)$ and $\ell_\infty$. By Corollary~\ref{c:6.2}, the Banach ring $C(K,\IF)$ is Hirsch and so is its isomorphic copy $\ell_\infty$. Since the Hirschness is inherited by subrings, the subrings $c$ and $c_0$ of the Hirsch ring $\ell_\infty$ are Hirsch.
\end{proof}

Observe that for every real number $p\in[1,\infty)$ the Banach space $\ell_p$ endowed with the pointwise multiplication of functions is a Banach ring. Theorem~\ref{t:lpq} implies

\begin{corollary}\label{c:12} For every $p\in[1,2]$ the Banach ring $\ell_p$ is Hirsch.
\end{corollary}


\begin{problem}\label{prob2} Is the Banach ring $\ell_p$ Hirsch for some/any $p\in(1,\infty)$?
\end{problem}

By a {\em Hilbert ring} we understand a Hilbert space $H$ endowed with a continuous associative operation $\cdot :H\times H\to H$ turning $H$ into a topological ring. An example of a Hilbert ring is the Banach ring $\ell_2$, which is Hirsch by Corollary~\ref{c:12}.

\begin{problem} Is every (commutative) Hilbert ring Hirsch?
\end{problem}

\subsection{Non-Hirsch Banach rings of operators}

Let $p\in [1,\infty]$ and $L(\ell_p)$ be the Banach space of bounded linear operators on $\ell_p$,
endowed with the operator norm $\|a\|=\sup\{\|a(x)\|:x\in\ell_p$ with $\|x\|_{\ell_p}\le 1\}$.
Let $L^{\!\circ}(\ell_p)$ (resp. $L^{\!\bullet}(\ell_p)$) be the Banach space $L(\ell_p)$ endowed with the multiplication $\circ$ (resp. $\bullet$) defined by $(a\circ b)(x)=a(b(x))$ (resp. $(a\bullet b)(x)=b(a(x))$) for any $a,b\in L(\ell_p)$ and $x\in\ell_p$. It is easy to check that both $L^{\!\circ}(\ell_p)$ and $L^{\!\bullet}(\ell_p)$ are Banach rings.

\begin{proposition}\label{p1} For any $p\in [1,\infty)$, the Banach ring $L^{\!\bullet}(\ell_p)$ is not Hirsch.
\end{proposition}

\begin{proof} Let $(e_i)_{i\in\w}$ be the standard basis of the Banach space $\ell_p$. For every $i\in\w$ let $a_i:\ell_p\to\ell_p$ be the operator assigning to each $(x_n)_{n\in\w}\in\ell_p$ the vector $\frac{x_i}{\sqrt[p]{i+1}}e_i$. Since the sequence $\big(\frac1{\sqrt[p]{i+1}}\big)_{i\in\w}$ converges to zero,
the series $\sum_{i\in\w} a_i$ is unconditionally Cauchy in $L(\ell_p)$, so it unconditionally converges
by Proposition~\ref{CauchyToUnCon}.
Now let $U=\{a\in L(\ell_p):\|a\|\le 1\}$ and $V$ be any open neighborhood of $0$ in $L(\ell_p)$.
Find  $\varepsilon>0$ such that $\e U\subseteq V$. Since the harmonic series $\sum_{i\in\w}
\frac1{i+1}$ diverges, there exists a positive integer $N$ such that
$\varepsilon^p\sum_{i=0}^N\frac1{i+1}>1$. For each $i\le N$, let $f_i:\ell_p\to\ell_p$ be the
linear operator assigning to each $(x_n)_{n\in\w}\in\ell_p$ the vector $\e x_0 e_i$. Observe that
$\{f_i\}_{i=0}^N\subseteq \e U\subseteq V$ but $$\Big\|\sum_{i=0}^N f_i{\bullet}
a_i\Big\|^p\ge\Big\|\Big(\sum_{i=0}^N f_i{\bullet} a_i\Big)(e_0) \Big\|^p_{\ell_p}=\varepsilon^p\sum_{i=0}^N
\Big(\tfrac1{\sqrt[p]{i+1}}\Big)^p>1.$$ So $\sum_{i=0}^N f_i{\bullet}a_i\not\in U$ and the Banach
ring $L^{\!\bullet}(\ell_p)$ is not Hirsch.
\end{proof}

\begin{proposition}\label{p1'} For any $p\in (1,\infty]$, the ring $L^{\!\circ}(\ell_p)$ is not Hirsch.
\end{proposition}
\begin{proof} Let $(e_i)_{i\in\w}$ be the standard basis of the Banach space $\ell_p$. For every $i\in\w$, let $a_i:\ell_p\to\ell_p$ be the operator assigning to each sequence $x=(x_n)_{n\in\w}\in\ell_p$ the vector $a_i(x)=\frac{x_i}{\ln (i+2)}e_i$. Since the sequence $\big(\frac1{\ln (i+2)}\big)_{i\in\w}$ converges to zero,
the series $\sum_{i\in\w} a_i$ is unconditionally Cauchy, so it unconditionally converges by Proposition~\ref{CauchyToUnCon}.

Let $U=\{f\in L(\ell_p):\|f\|\le 1\}$ be the closed unit ball in the Banach space $L(\ell_p)$. Given any neighborhood of zero $V\subseteq L(\ell_p)$, find $\varepsilon>0$ such that $\e U\subseteq V$.

 Let $e=\sum_{n\in\w}\tfrac1{n+2}e_n$.
Since $p>1$, $e\in\ell_p$, see, for instance,~\cite[367.5]{Fich}.

Since the series $\sum_{i=2}^\infty \frac1{i\ln i}$ diverges (see, for
instance,~\cite[367.6]{Fich}), there exists an integer number $N$ such that
$\varepsilon\sum_{i=0}^N \tfrac1{(i+2)\ln(i+2)}>\|e\|_{\ell_p}$. For each $i\le N$ let $f_i:\ell_p\to\ell_p$
be the operator assigning to each vector $x=(x_n)\in\ell_p$ the vector $f_i(x)=\e x_ie_0$.  Then
$\{f_i\}_{i=0}^N\subseteq \e U\subseteq V$ and $$\Big\|\sum_{i\in\w} f_i{\circ}a_i
\Big\|\cdot\big\|e\big\|_{\ell_p}\ge\Big\|\Big(\sum_{i\in\w} f_i{\circ}a_i\Big)(e)
\Big\|_{\ell_p}=\varepsilon\sum_{i=0}^N \tfrac1{(i+2)\ln(i+2)}
>\|e\|_{\ell_p}.$$ So $\sum_{i\in\w} f_i{\circ}a_i\not\in U$ and hence $L^{\!\circ}(\ell_p)$ is not Hirsch.
\end{proof}

\begin{question} Are the Banach rings $L^{\!\bullet}(\ell_\infty)$ and $L^{\!\circ}(\ell_1)$ Hirsch?
\end{question}

\subsection{Non-Hirsch commutative Banach rings}

In this subsection we present a simple example of a commutative Banach ring, which is not Hirsch. For a real number $p\in[1,\infty)$ we denote by $\ell_p$  the real Banach space of functions $x:\w\to\IR$ with finite norm
$$\|x\|_{\ell_p}=\Big(\sum_{k\in\w}|x(k)|^p\Big)^{\frac1p}<\infty.$$ Let $\Re:\IC\to\IR$ and $\Im:\IC\to\IR$ be functions of taking the real and imaginary parts of a complex number, respectively.

Given two real numbers $p,q\in[1,\infty)$, consider the subring
$$\ell_p\oplus i\ell_q=\{z\in\IC^\w:\Re\circ z\in\ell_p,\;\;\Im\circ z\in\ell_q\}$$of the commutative ring $\IC^\w$ of all complex-valued sequences. The ring $\ell_p \oplus  i\ell_q$ is a Banach ring with respect to the norm
$$\|z\|=\max\{\|\Re\circ z\|_p,\|\Im\circ z\|_q\}.$$

\begin{theorem}\label{t:pq}  For any distinct real numbers $p,q\in[1,\infty)$, the commutative Banach ring $\ell_p \oplus i\ell_q$ is not Hirsch.
\end{theorem}

\begin{proof} For every $n\in\w$ let $e_n:\w\to\{0,1\}$ be the unique function such that $e^{-1}_n(1)=\{n\}$.

Let $m=\min\{p,q\}$ and $M=\max\{p,q\}$. Let
$$b=\begin{cases}1&\mbox{if $p<q$};\\
0&\mbox{if $q<p$}.
\end{cases}
$$
Then $i^b=i$ if $p<q$ and $i^b=1$ if $q<p$.

Observe that the series $\sum_{k\in\w}(k+1)^{-1/m}e_k$ converges unconditionally in the Banach space $\ell_M$ but diverges in $\ell_m$. The choice of $b$ ensures that the series $$\sum_{k\in\w}x_k=\sum_{k\in\w}i^b(k+1)^{-1/m}e_k$$ converges unconditionally in the Banach space $\ell_p\oplus i\ell_q$.

Assuming for a contradiction that the Banach ring $\ell_p\oplus i\ell_q$ is Hirsch, for the neighborhood $U=\{z\in\ell_p\oplus i\ell_q:\|z\|<1\}$  of zero, find a neighborhood $V\subseteq\ell_p\oplus i\ell_q$ of zero such that for every $n\in\w$ and elements $(a_k)_{k\in n}\in V^n$ we have $\sum_{k\in n}a_kx_k\in U$.
Find $\e>0$ such that $\{z\in\ell_p \oplus i\ell_q:\|z\|\le\e\}\subseteq V$. Since the series $\sum_{k\in\w}(k+1)^{-1}$ diverges, there exists $n\in\w$ such that $$\Big(\sum_{k\in n}(k+1)^{-1}\Big)^{\frac1m}>\frac1\e.$$

For every $k\in n$ let $a_k=\e ie_k\in V$. Then $$
\begin{aligned}
\big\|\sum_{k\in n}a_kx_k\big\|&=
\big\|\sum_{k\in n}(\e i e_k)(i^b(k+1)^{-1/m}e_k)\big\|=\big\|\sum_{k\in n}\e i^{b+1}(k+1)^{-1/m}e_k\big\|=\\
&=\big\|\sum_{k\in n}\e (k+1)^{-1/m}e_k\big\|_{\ell_m}=\e\Big(\sum_{k\in n}(k+1)^{-1}\Big)^{\frac1m}>1
\end{aligned}
$$
and hence $\sum_{k\in n}a_kx_k\notin U$, which contradicts the choice of the neighborhood $V$.
\end{proof}

\section{Acknowledgements}

The authors thank to Dylan Hirsch for asking an interesting question on {\tt Mathematics Stack Exchange}~\cite{Dch} and to Davide Giraudo for the help in the attempts to solve Problem~\ref{prob2}.


\begin{thebibliography}{}

\bibitem{ArhTka}
A.~Arhangel'ski\u{\i}, M.~Tkachenko,
{\it Topological Groups and Related Structures},
Atlantis Press, 2008.

\bibitem{BanKad} T.~Banakh, V.~Kadets, {\em Multiplications preserving unconditional convergence in Banach spaces}, Axioms (to appear); \url{arxiv.org/abs/2111.14253}.

\bibitem{BL} T.~Banakh, A.~Leiderman, {\em The strong Pytkeev property in topological spaces}, Topology Appl. {\bf 227} (2017) 10--29.





\bibitem{Fern} C.S.~Fernandez, {\em The closed graph theorem for multilinear mappings}, Internat. J. Math. Math. Sci. {\bf 19}:2 (1996), 407--408.


\bibitem{Fich}
G. Fichtenholz,
{\it Differential and integral calculus},
v. II, 7th ed., M.:Nauka, 1970, in Russian.

\bibitem{Dch}
D. Hirsch ({\tt dch}), {\it Unconditional convergence of a sum of elements in a complete Hausdorff topological ring}, 
\url{math.stackexchange.com/q/3702472}.

\bibitem{HM}
K.H. Hofmann, S.A. Morris,
{\it The Structure of Compact Groups},
3rd ed., de Gruyter, Berlin, 2013.


\bibitem{NV} I.~Netuka, Ji\v r\'\i\ Vesel\'y, {\em An inequality for finite sums in ${\bf R}\sp{m}$}, \v Casopis P\v st. Mat. {\bf 103}:1 (1978), 73--77.

\bibitem{Pon}
L.S. Pontrjagin,
{\it Continuous groups},
2nd ed., M., 1954, in Russian.
\end{thebibliography}
\end{document}